\documentclass{article}%
\usepackage{amssymb}
\usepackage{amsmath}
\usepackage{amsfonts}
\usepackage{graphicx}%
\providecommand{\U}[1]{\protect \rule{.1in}{.1in}}
\newtheorem{theorem}{Theorem}

\newtheorem{proposition}[theorem]{Proposition}

\begin{document}

\title{\textbf{Surfaces of revolution}\\ \textbf{satisfying }$\triangle^{III}\boldsymbol{x}=A\boldsymbol{x}$}
\author{\textbf{Stylianos Stamatakis, Hassan Al-Zoubi\medskip}\\ \emph{Department of Mathematics, Aristotle University of Thessaloniki}\\ \emph{GR-54124 Thessaloniki, Greece}\\ \emph{e-mail: stamata@math.auth.gr}}
\date{}
\maketitle

\begin{abstract}
\noindent We consider surfaces of revolution in the three-dimensional
Euclidean space which are of coordinate finite type with respect to the third
fundamental form $III$, i.e., their position vector $\boldsymbol{x}$ satisfies
the relation \textbf{ }$\triangle^{III}\boldsymbol{x}=A\boldsymbol{x}%
$,\textbf{ }where $A$ is a square matrix of order 3. We show that a surface of
revolution satisfying the preceding relation is a catenoid or part of a sphere.



\end{abstract}

\section{Introduction}

\noindent Let $\boldsymbol{x}=\boldsymbol{x}(u^{1},u^{2})$ be a regular
parametric representation of a surface $S$ in the Euclidean space $%
\mathbb{R}
^{3}$ which does not contain parabolic points. For two sufficient
differentiable functions $f(u^{1},u^{2})$ and $g(u^{1},u^{2})$ the first
Beltrami operator with respect to the third fundamental form $III=e_{ij}%
du^{i}du^{j}$ of $S$ is defined by%
\[
\nabla^{III}(f,g)=e^{ij}f_{/i}g_{/j},
\]
where $f_{/i}:=\frac{\partial f}{\partial u^{i}}$ and $e^{ij}$ denote the
components of the inverse tensor of $e_{ij}$. The second Beltrami differential
operator with respect to $III$ is defined by \footnote{with sign convention
such that $\triangle=-\frac{\partial^{2}}{\partial x^{2}}-\frac{\partial^{2}%
}{\partial y^{2}}$ for the metric $ds^{2}=dx^{2}+dy^{2}$}
\begin{equation}
\triangle^{III}f=\frac{-1}{\sqrt{e}}\left(  \sqrt{e}e^{ij}f_{/i}\right)
_{/j}\label{1}%
\end{equation}
($e:=\det(e_{ij})).$ In \cite{Stamatakis} we showed the relation%
\begin{equation}
\triangle^{III}\boldsymbol{x}=\nabla^{III}(\frac{2H}{K},\boldsymbol{n}%
)-\frac{2H}{K}\boldsymbol{n},\label{2}%
\end{equation}
where $\boldsymbol{n}$ is the unit normal vectorfield, $H$ the mean curvature
and $K$ the Gaussian curvature of $S$. Moreover we proved that a surface
satisfying the condition%
\[
\triangle^{III}\boldsymbol{x}=\lambda \boldsymbol{x},\quad \lambda \in%
\mathbb{R}
,
\]
i.e., a surface $S:\boldsymbol{x}=\boldsymbol{x}(u^{1},u^{2})$ \textit{for
which all coordinate functions are eigenfunctions} of $\triangle^{III}$
\textit{with the same eigenvalue} $\lambda$, is part of a sphere ($\lambda=2$)
or a minimal surface ($\lambda=0$). Using terms of B.-Y. Chen's theory of
finite type surfaces \cite{Chen} the above result can be expressed as follows:
\textit{A surface} $S$ \textit{in} $%
\mathbb{R}
^{3}$ \textit{is of }$III$-\textit{type} $1$ (\textit{or of null}
$III$-\textit{type} $1$) \textit{if and only if} $S$ \textit{is part of a
sphere }(\textit{or a minimal surface})\textit{.}

In general a surface $S$ is said to be \textit{of finite type }with respect to
the fundamental form $III$ or, briefly, \textit{of finite }$III$%
\textit{-type}, if the position vector $\boldsymbol{x}$ of $S$ can be written
as a finite sum of nonconstant eigenvectors of the operator $\triangle^{III}$,
that is if%
\begin{equation}
\boldsymbol{x}=\boldsymbol{c}+\boldsymbol{x}_{1}+\boldsymbol{x}_{2}%
+\ldots+\boldsymbol{x}_{m},\quad \triangle^{III}\boldsymbol{x}_{i}=\lambda
_{i}\boldsymbol{x}_{i},\quad i=1,\ldots,m, \label{3}%
\end{equation}
where $\boldsymbol{c}$ is a constant vector and $\lambda_{1},\ldots
,\lambda_{m}$ are eigenvalues of $\triangle^{III}$. When there are exactly $k$
nonconstant eigenvectors $\boldsymbol{x}_{1},\ldots,\boldsymbol{x}_{k}$
appearing in (\ref{3}) which all belong to different eigenvalues $\lambda
_{1},\ldots,\lambda_{k}$, then $S$ is said to be of $III$-\textit{type} $k$;
when $\lambda_{i}=0$ for some $i=1,\ldots,k$, then $S$ is said to be of
\textit{null }$III$-\textit{type} $k$.

The only known surfaces of finite $III$-type are parts of spheres, the minimal
surfaces and the parallel of the minimal surfaces (which are actually of null
$III$-type $2$, see \cite{Stamatakis}).

In this paper we want to determine the connected surfaces of revolution $S$ in
$%
\mathbb{R}
^{3}$ which are \textit{of coordinate finite }$III$-\textit{type}, i.e., their
position vectorfield $\boldsymbol{x}(u^{1},u^{2})$ satisfies the condition%
\begin{equation}
\triangle^{III}\boldsymbol{x}=A\boldsymbol{x},\quad A\in M(3,3),\label{4}%
\end{equation}
where $M(m,n)$ denotes the set of all matrices of the type ($m,n$).

Coordinate finite type surfaces with respect to the first fundamental form $I$
were studied in \cite{Dillen} and \cite{Garay}. In the last paper O. Garay
showed that the only complete surfaces of revolution in $%
\mathbb{R}
^{3}$, whose component functions are eigenfunctions of their Laplacian are the
catenoids, the spheres and the circular cylinders, while F. Dillen, J. Pas and
L. Verstraelen proved in \cite{Dillen} that the only surfaces in $%
\mathbb{R}
^{3}$ satisfying%

\[
\triangle^{I}\boldsymbol{x}=A\boldsymbol{x}+B,\quad A\in M(3,3),\quad B\in
M(3,1),
\]
are the minimal surfaces, the spheres and the circular cylinders.

Our main result is the following

\begin{proposition}
A surface of revolution $S$ satisfies (\ref{4}) if and only if $S$ is a
catenoid or part of a sphere.
\end{proposition}

We first show that the mentioned surfaces indeed satisfy the condition
(\ref{4}).

A. On a catenoid the mean curvature vanishes, so, by virtue of (\ref{2}),
$\triangle^{III}\boldsymbol{x}=0$. Therefore a catenoid satisfies (\ref{4}),
where $A$ is the null matrix in $M(3,3)$.

B. Let $S$ be part of a sphere of radius $r$ centered at the origin. Then%
\[
H=\frac{1}{r},\quad K=\frac{1}{r^{2}},\quad \boldsymbol{n}=-\frac{1}%
{r}\boldsymbol{x}.
\]
So, by (\ref{2}), it is $\triangle^{III}\boldsymbol{x}=2\boldsymbol{x}$.
Therefore $S$ satisfies (\ref{4}) whith $A=2I_{3},$ where $I_{3}$ is the
identity matrix in $M(3,3)$.

\section{Proof of the main theorem}

\noindent Let $C$ be the profile curve of a surface of revolution $S$ of the
differentiation class $C^{3}$. We suppose that (a) $C$ lies on the
$(x_{1},x_{3})$-plane, (b) the axis of revolution of $S$ is the $x_{3}$-axis
and (c) $C$ is parametrized by its arclength $s$. Then $C$ admits the
parametric representation%
\[
\boldsymbol{r}(s)=(f(s),0,g(s)),\quad s\in J
\]
($J\subset%
\mathbb{R}
$ open interval), where $f(s),g(s)\in C^{3}(J)$. The position vector of $S$ is
given by%
\[
\boldsymbol{x}(s,\theta)=(f(s)\cos \theta,f(s)\sin \theta,g(s)),\quad s\in
J,\quad \theta \in \lbrack0,2\pi).
\]
Putting $f(s)%
\acute{}%
:=\frac{df(s)}{ds}$ we have because of (c)
\begin{equation}
f~%
\acute{}%
~^{2}+g~%
\acute{}%
~^{2}=1\quad \forall s\in J.\label{5}%
\end{equation}
Furthermore it is $f~%
\acute{}%
\cdot g~%
\acute{}%
\neq0,$ because otherwise $f=const.$ or $g=const.$ and $S$ would be a circular
cylinder or part of a plane, respectively. Hence $S$ would consist only of
parabolic points, which has been excluded. In view of (\ref{5}) we can put%
\begin{equation}
f~%
\acute{}%
=\cos \varphi,\quad g~%
\acute{}%
=\sin \varphi,\label{6}%
\end{equation}
where $\varphi$ is a function of $s$. Then the unit normal vector of $S$ is
given by%
\[
\boldsymbol{n}=(-\sin \varphi \cos \theta,~-\sin \varphi \sin \theta,~\cos \varphi).
\]
The components $h_{ij}$ and $e_{ij}$ of the the second and the third
fundamental tensors in (local) coordinates are the following%
\[
h_{11}=\varphi%
\acute{}%
,\quad h_{12}=0,\quad h_{22}=f\sin \varphi,
\]%
\begin{equation}
e_{11}=\varphi%
\acute{}%
~^{2},\quad e_{12}=0,\quad e_{22}=\sin^{2}\varphi,\label{7}%
\end{equation}
hence \cite{Huck}%

\begin{equation}
\frac{2H}{K}=h_{ij}e^{ij}=\frac{1}{\varphi%
\acute{}%
}+\frac{f}{\sin \varphi}.\label{8}%
\end{equation}
From (\ref{1}) and (\ref{7}) we find for a sufficient differentiable function
$u=u(s,\theta)$ defined on $J\times \lbrack2\pi,0)$
\begin{equation}
\triangle^{III}u=-\frac{u\text{%
\'{}%
\'{}%
}}{\varphi%
\acute{}%
^{~2}}+\left(  \frac{\varphi \text{%
\'{}%
\'{}%
}}{\varphi%
\acute{}%
^{~2}}-\frac{\cos \varphi}{\sin \varphi}\right)  \frac{u%
\acute{}%
}{\varphi%
\acute{}%
}-\frac{u_{/\theta \theta}}{\sin^{2}\varphi}.\label{9}%
\end{equation}
Consindering the following functions of $s$%
\begin{equation}
P_{1}=R\sin \varphi-\frac{\cos \varphi}{\varphi%
\acute{}%
}R~%
\acute{}%
,\quad P_{2}=-R\cos \varphi-\frac{\sin \varphi}{\varphi%
\acute{}%
}R~%
\acute{}%
,\label{11}%
\end{equation}
where we have put for simplicity $R:=\frac{2H}{K}$, and applying (\ref{9}) on
the coordinate functions $x_{i}$, $i=1,2,3,$ of the position vector
$\boldsymbol{x}$ we find%
\begin{equation}
\triangle^{III}x_{1}=P_{1}\cos \theta,\quad \triangle^{III}x_{2}=P_{1}\sin
\theta,\quad \triangle^{III}x_{3}=P_{2}.\label{10}%
\end{equation}
So we have: 

(a) \textit{The coordinate functions }$x_{1},x_{2}$ \textit{are both
eigenfunctions of} $\triangle^{III}$ \textit{belonging to the same eigenvalue
if and only if for some real constant} $\lambda$ \textit{holds}%

\[
\lambda f=R\sin \varphi-\frac{\cos \varphi}{\varphi%
\acute{}%
}R~%
\acute{}%
.
\]

(b) \textit{The coordinate function }$x_{3}$ \textit{is an eigenfunction of}
$\triangle^{III}$ \textit{if and only if for some real constant} $\mu$
\textit{holds}%
\[
\mu g=-R\cos \varphi-\frac{\sin \varphi}{\varphi%
\acute{}%
}R~%
\acute{}%
.
\]
We denote by $a_{ij},i,j=1,2,3,$ the entries of the matrix $A.$ By using
(\ref{10}) condition (\ref{4}) is found to be equivalent to the following
system%
\begin{equation}
\left \{
\begin{array}
[c]{c}%
P_{1}\cos \theta=a_{11}f\cos \theta+a_{12}f\sin \theta+a_{13}~g\\
P_{1}\sin \theta=a_{21}f\cos \theta+a_{22}f\sin \theta+a_{23}~g\\
\quad \quad P_{2}=a_{31}f\cos \theta+a_{32}f\sin \theta+a_{33}~g
\end{array}
\right.  .\label{12}%
\end{equation}
Since $\sin \theta,\cos \theta$ and 1 are linearly independent functions of
$\theta,$ we obtain from (\ref{12}$_{3}$) $a_{31}=a_{32}=0.$ On
differentiating (\ref{12}$_{1}$) and (\ref{12}$_{2}$) twice with respect to
$\theta$ we have%
\[
\left \{
\begin{array}
[c]{c}%
P_{1}\cos \theta=a_{11}f\cos \theta+a_{12}f\sin \theta \\
P_{1}\sin \theta=a_{21}f\cos \theta+a_{22}f\sin \theta
\end{array}
\right.  .
\]
Thus $a_{13}g=a_{23}g=0,$ so that $a_{13}$ and $a_{23}$ vanish. The system
(\ref{12}) is equivalent to the following%
\[
\left \{
\begin{array}
[c]{c}%
\left(  P_{1}-a_{11}f\right)  \cos \theta-a_{12}f\sin \theta=0\\
(P_{1}-a_{22}f)\sin \theta-a_{21}f\cos \theta=0\\
\quad \quad \quad \quad \quad \quad \quad \quad P_{2}-a_{33}g=0
\end{array}
\right.  .
\]
But $\sin \theta$ and $\cos \theta$ are linearly independent functions of
$\theta$, so we finally obtain $a_{12}=a_{21}=0,a_{11}=a_{22}$ and
$P_{1}=a_{11}f.$ Putting $a_{11}=a_{22}=\lambda$ and $a_{33}=\mu$ we see that
the system (\ref{12}) reduces now to the following equations%
\begin{equation}
P_{1}=\lambda f,\quad P_{2}=\mu g.\label{13}%
\end{equation}
On account of (\ref{11}) and (\ref{13}) we are left with the system%
\begin{equation}
\left \{
\begin{array}
[c]{c}%
R=\lambda f\sin \varphi-\mu g\cos \varphi \quad \quad \\
R~%
\acute{}%
=-\varphi~%
\acute{}%
(\lambda f\cos \varphi+\mu g\sin \varphi)
\end{array}
\right.  .\label{14}%
\end{equation}
On differentiating (\ref{14}$_{1}$) with respect to $s$ we find, by virtue of
(\ref{6}),%
\begin{equation}
R~%
\acute{}%
=\frac{\lambda-\mu}{2}\sin \varphi \cos \varphi.\label{15}%
\end{equation}
We distinguish the following cases:

\textit{Case I.} Let $\lambda=\mu$.

Then (\ref{15}) reduces to $R~%
\acute{}%
=0.$

\textit{Subcase Ia}. Let $\lambda=\mu=0$. From (\ref{14}$_{1}$) we obtain
$R=0$, i.e., $H=0$. Consequently $S$, being a minimal surface of revolution,
is a catenoid.

\textit{Subcase Ib}. Let $\lambda=\mu \neq0$.

Then from (\ref{6}), (\ref{14}$_{2}$) and $R~%
\acute{}%
=0$ we have $f\cdot f~%
\acute{}%
+g\cdot g~%
\acute{}%
=0,$ i.e., $\left(  f^{2}+g^{2}\right)
\acute{}%
=0.$ Therefore $f^{2}+g^{2}=const.$ and $S$ is obviously part of a sphere.

\textit{Case II.} Let $\lambda \neq \mu$. From (\ref{14}$_{2}$), (\ref{15}) we
find firstly%
\begin{equation}
\frac{1}{\varphi%
\acute{}%
}=\frac{2\left(  \lambda f\cos \varphi+\mu g\sin \varphi \right)  }{\left(
\mu-\lambda \right)  \sin \varphi \cos \varphi}. \label{16}%
\end{equation}
From this and (\ref{8}) we obtain%
\[
R=\frac{\lambda+\mu}{\left(  \mu-\lambda \right)  \sin \varphi}f+\frac{2\mu
}{\left(  \mu-\lambda \right)  \cos \varphi}g.
\]
Hence, by virtue of (\ref{14}$_{1}$),%
\begin{equation}
af+bg=0, \label{17}%
\end{equation}
where%
\begin{equation}
a=\lambda \sin \varphi+\frac{\lambda+\mu}{\left(  \lambda-\mu \right)
\sin \varphi},\quad b=\frac{2\mu}{\left(  \lambda-\mu \right)  \cos \varphi}%
-\mu \cos \varphi. \label{18}%
\end{equation}
We note that $\mu \neq0$, since for $\mu=0$ we have%
\[
a=\frac{\lambda \sin^{2}\varphi+1}{\sin \varphi},\quad b=0,
\]
and relation (\ref{17}) becomes%
\[
\frac{\lambda \sin^{2}\varphi+1}{\sin \varphi}f=0,
\]
whence it follows $\lambda \sin^{2}\varphi+1=0,$ a contradiction.

On differentiating (\ref{17}) with respect to $s$ and taking into account
(\ref{16}) we obtain%
\begin{equation}
a_{1}\frac{f}{\sin \varphi}+b_{1}\frac{g}{\cos \varphi}=0, \label{19}%
\end{equation}
where%
\begin{equation}
a_{1}=\lambda(\lambda-\mu)^{2}\sin^{4}\varphi+(\lambda-\mu)(\lambda \mu
-\lambda^{2}+3\lambda+\mu)\sin^{2}\varphi-(\lambda+\mu)(3\lambda-\mu),
\label{20}%
\end{equation}%
\begin{equation}
b_{1}=\mu \left[  \left(  \lambda-\mu \right)  ^{2}\sin^{4}\varphi+\left(
\lambda-\mu \right)  \left(  \mu-\lambda+4\right)  \sin^{2}\varphi-2\left(
\lambda+\mu \right)  \right]  . \label{21}%
\end{equation}
By eliminating now the functions $f$ and $g$ from (\ref{17}) and (\ref{19})
and taking into account (\ref{18}), (\ref{20}) and (\ref{21}) we find%
\[
\lambda(\lambda-\mu)^{2}\sin^{4}\varphi+(\lambda-\mu)(\lambda \mu-\lambda
^{2}+5\lambda+\mu-2)\sin^{2}\varphi+(\lambda+\mu)(\mu-3\lambda+4)=0.
\]
Consequently%
\[
\lambda \left(  \lambda-\mu \right)  ^{2}=0,~\left(  \lambda-\mu \right)  \left(
\lambda \mu-\lambda^{2}+5\lambda+\mu-2\right)  =0,~\left(  \lambda+\mu \right)
\left(  \mu-3\lambda+4\right)  =0.
\]
From the first equation we have $\lambda=0$. Then, the other two become as
follows%
\[
\mu-2=0,\quad \mu+4=0,
\]
which is a contradiction.

So the proof of the theorem is completed.

\end{document}